\documentclass[11pt]{amsart} 
\usepackage{pb-diagram,lamsarrow,pb-lams,epic,eepic}
\usepackage{amsmath,amsthm,amssymb,amsfonts,latexsym,fullpage}

\newtheorem{theorem}{Theorem}[section]

\newtheorem{proposition}[theorem]{Proposition}   

\theoremstyle{definition}
\newtheorem{definition}[theorem]{Definition}
\newtheorem{remark}[theorem]{Remark}
\newtheorem{example}[theorem]{Example}

\def\H{{\rm H}}
\def\n{{\mathbb N}}
\def\Tor{{\rm Tor}}
\def\ann{{\rm Ann}}
\def\im{{\rm Im}}
\def\pd{{\rm Pd}}
\def\depth{{\rm depth}}

\def\y{{\mathbf{y}}}  
\def\phir{{\phi^r}}
\def\rphir{{^{\phi^r}\!\! R}}
\def\fphir{{^{\phi^r}\!\! F_{\bullet}}}
\def\dphir{{^{\phi^r}\!\! d}}

\begin{document}

\title[Frobenius powers of non-complete intersections]
{Frobenius powers of non-complete
intersections}  
\author[M. R. Kantorovitz]{Miriam Ruth Kantorovitz}
\address{Dept. of Math, University of California,
Berkeley, CA 94720 }  
\email{ ruth@math.berkeley.edu}
\thanks{Partially supported by an NSF grant}
\maketitle

\section*{Introduction}
The purpose of this paper is to address a number of issues raised by 
Avramov and Miller in a recent paper~\cite{AM}.

Let $(R,m,k)$ be a Noetherian local ring of characteristic $p>0$ with
residue field $k$, and let $\phi : R \to R$ be the
the Frobenius homomorphism defined by $\phi(a)=a^p$.
For $r \geq 1$, we denote by $\rphir$ the $R$-module structure on $R$ via
$\phir$. That is, for $a\in R$ and $b \in \rphir$, $a\cdot b =
a^{p^r}b$. When $R$ is a regular ring, $\rphir$ is flat; in
fact, this condition characterizes regular rings~\cite{Ku}.  
When $R$ is a complete intersection,
Avramov and Miller~\cite{AM} proved that 
$\Tor^R_*(-,\rphir)$ is rigid in the following sense.
\begin{theorem}\label{AM1}(cf. \cite[main theorem]{AM})
If $R$ is a complete intersection and $M$ is an $R$-module such that 
$\Tor^R_j(M,\rphir)=0$ holds for some $j,r \geq 1$, then 
$\Tor^R_n(M,\rphir)=0$ for all $n \geq j$. Furthermore, if $M$ is
finitely generated then $M$ has finite projective dimension.
\end{theorem}
This theorem shows similarity between the functors $\Tor^R_*(-,\rphir)$ 
and $\Tor^R_*(-,k)$ in terms of rigidity.
If $M$ is an $R$-module of finite length, $\ell_R(M)$, then, in addition,
the following relationship holds between the lengths of the homology modules
$\Tor^R_*(M,\rphir)$ 
and $\Tor^R_*(M,k)$.
\begin{theorem}\label{AM2}(cf. \cite[main theorem]{AM})
If $R$ is a complete intersection and $M$ is an $R$-module of finite
length and infinite projective dimension then for each $r \geq 1$,
both
$$
\lim_{s \to \infty}\frac{
\ell_R(\Tor^R_{2s}(M,\rphir))}{
\ell_R(\Tor^R_{2s}(M,k))}
\quad \mbox{ and } \quad
\lim_{s \to \infty}\frac{
\ell_R(\Tor^R_{2s+1}(M,\rphir))}{
\ell_R(\Tor^R_{2s+1}(M,k))}
$$ 
are rational numbers, and at least one of them is positive.
\end{theorem}
Avramov and Miller wondered whether the complete
intersection assumption on the ring $R$ is necessary in
Theorems~\ref{AM1} and~\ref{AM2}. In this paper we give an
answer to their question. 
First we give a class of depth zero rings 
(which include non-complete intersections), for which
 rigidity of $\Tor^R_*(-,\rphir)$ holds.
The simplest examples of rings 
 for which the conclusions of
Theorems~\ref{AM1}
and~\ref{AM2} hold are Artinian rings $(R,m)$ with $m^p=0$.
A simple argument is used to show that if $R$ satisfies the
condition 
$$
(1) \qquad (0:m^p)_R \not\subseteq m^p,
$$
then strong rigidity of $\Tor$ holds for finitely
generated $R$-modules. That is, if $\Tor_j^R(M, \rphir)=0$ then $M$ is
projective. 
However, when $\depth (R) >0$, we show that 
$\Tor_j^R(-, \rphir)$ is not rigid  (non-vacuously) in general 
(see Proposition~\ref{prop-not-rigid}) and hence the complete 
intersection assumption in Theorem~\ref{AM1} is necessary.
We conclude the paper with a few examples of non-complete
intersections which satisfy condition $(1)$.
\subsection*{Acknowledgments} I thank David Eisenbud for a number of
useful conversations and Claudia Miller for her comments.
\section{Depth zero rings}
Let $M$ be a finitely generated $R$-module. 
In order to compute the
homology modules $\Tor^R_*(M,\rphir)$, 
 choose a minimal free resolution of $M$,  
$$
F_{\bullet}: \qquad \dots \to F_{n+1} \stackrel{d}{\to} F_{n}   \to
\dots ,
$$
where
 $F_n = R^{\oplus l_n}$ for some 
$l_n \in \n$, and the differential maps, $d = (d_{ij})$, are
represented (with respect to the standard set of generators) by matrices
with entries in $m$.
Applying the functor $- \otimes_R \rphir$ to the resolution $F_{\bullet}$, we
obtain the complex 
$$
F_{\bullet}\otimes_R \rphir : \qquad \dots \to 
F_{n+1} \otimes_R \rphir 
\stackrel{d \otimes 1}{\longrightarrow} 
F_n \otimes_R \rphir  \to  \dots ,
$$
 the homology of which computes $\Tor^R_*(M,\rphir)$.
After making the standard identifications, this complex is homologically
equivalent to the complex
$$
^{\phi^r}\!\!F_{\bullet}: \qquad \dots \to R^{\oplus l_{n+1}} 
\stackrel{\dphir = (d_{ij}^{p^r})}{\longrightarrow} 
R^{\oplus l_n}   \to  \dots .
$$
That is, $\Tor^R_*(M,\rphir)=H_n(^{\phi^r}\!\!F_{\bullet})$, where 
$^{\phi^r}\!\!F_{\bullet}$ is the complex obtained from $F_{\bullet}$
by raising the entries in the differential $d=(d_{ij})$ to the
$p^r$-th power.
Note that since $d_{ij} \in m$, the image of the $n$-th differential, $\im
(\dphir)$, is contained in $m^{p^r}F_n$.
\begin{proposition} \label{prop1}
Let $(R,m)$ be a local ring of characteristic $p$ satisfying the
following condition
\begin{equation} \label{e:condition}
(0:m^p)_R := \{ x\in R | m^p  \subseteq \ann (x)\} \not\subseteq m^p.
\end{equation}
If $M$ is a finitely generated $R$-module such that 
$\Tor^R_j(M,\rphir)=0$ for some $j,r \geq 1$ then $M$ is projective.
\end{proposition}
\begin{proof}
Using the
above notation, let $F_{\bullet} \to M$ be a minimal free resolution of $M$
with $F_j = R^{\oplus l_j}$ for some 
$l_j \in \n$. First we show that $F_j=0$. 
By~\eqref{e:condition}, there exist an element 
$w \in m \smallsetminus m^p$ such that $m^p \subset \ann(w)$. 
If $F_j \neq 0$, then in the
complex $(^{\phi^r}\!\!F_{\bullet}, \dphir)$, the element 
$\underline{w}=(w, w, \dots , w) \in F_j \smallsetminus m^{p^r}F_j$ satisfies 
$\dphir(\underline{w})=0$ since $(\dphir)_{ij} \in m^{p^r} \subset
\ann(w)$.
 However,
 $\im(\dphir) \subset m^{p^r}F_j$, and hence $\underline{w}$ is a
$j$-cycle which is not a boundary. Therefore $\Tor^R_j(M,\rphir) \neq
0$, which contradicts the assumption. Hence $F_j=0$, and the
projective dimension of $M$, $\pd (M)$, is finite.

Observe that condition~\eqref{e:condition} implies that the depth of
$R$ is zero. Hence by the Auslander-Buchsbaum formula, 
$\pd (M) + \depth (M) = \depth (R)$,
 $M$ is projective.
\end{proof}
\begin{remark}\label{remark0}
What we actually proved above is the following:
If $(R,m)$ satisfies~\eqref{e:condition} 
and $L_{\bullet}$ is a complex of finitely generated 
free $R$-modules with differentials having entries in $m$, then
$$
\H_j(L_{\bullet} \otimes_R \rphir) =0 \Leftrightarrow L_j=0.
$$
\end{remark}
Following~\cite{AM}, we define the integer $c_\y(R)$, associated to 
the ring R and to a maximal $R$-regular sequence $\y$.
\begin{definition}\label{def1}
When $\depth (R) =0$, define $c(R)$ to be the smallest integer $s$ such that 
$(0:m)_R = \{ x \in R| m \subset \ann (x) \} \not\subset m^s$.
Note that such an $s$ exists by Krull's intersection theorem
since $(0:m)_R \neq 0$.
When $\depth (R)=d >0$, let $\y = \{y_1, \ldots, y_d\}$ 
be a maximal $R$-regular sequence. Write 
 $\bar{R}$ for $\frac{R}{\y R}$ and $\bar{m}$ for 
the corresponding maximal ideal of $\bar{R}$.
Define $c_\y(R)$ to be the least positive integer $s$
such that $(0 : \bar{m})_{\bar{R}} \not \subset \bar{m}^s$.
\end{definition}
The same argument as in the proof of Proposition~\ref{prop1}
 gives the following rigidity result
for all local rings of depth zero. 
\begin{proposition}\label{prop-depth0}
Let $(R,m)$ be a local ring of depth zero and let 
$L_{\bullet}$ be a complex of finitely generated free $R$-modules with 
differentials having entries in $m$. Then for $r > \log_p c(R)$,
$$
\H_j(L_{\bullet} \otimes_R \rphir) =0 \Leftrightarrow L_j=0.
$$
Consequently, if $M$ is a finitely generated $R$-module such that 
$\Tor^R_j(M,\rphir)=0$ for some $j \geq 1$ and $r > \log_pc(R)$, 
then $M$ is projective.
\end{proposition}
For Artinian rings, 
Proposition~\ref{prop1} can be trivially extended to 
non-finitely generated modules since in that case, we still
have a notion of a minimal free resolutions. 
\begin{proposition}\label{prop2}
Let $(R,m)$ be an Artinian local ring of characteristic $p$ satisfying
$m^p=0$. If $M$ is an $R$-module such that $\Tor^R_j(M,\rphir)=0$ for
some $j,r \geq 1$ then $\pd (M) < j$. 
\end{proposition}
\begin{proof}
Let $F_{\bullet} \to M$ be a minimal free resolution of $M$. Then
in the complex $\fphir$, all the differentials are zero and 
therefore, $F_j=0$. 
\end{proof}

We now turn our attention to the relation between the lengths of the
homology modules
$\Tor^R_*(M,\rphir)$ and $\Tor^R_*(M,k)$ when $M$ is of finite length.
\begin{proposition}\label{prop3}
Let $(R,m)$ be an Artinian local ring of characteristic $p$ satisfying
$m^p=0$. Let $M$ be an $R$-module of finite
length and infinite projective dimension then for each $r \geq 1$,
$$
\lim_{s \to \infty}\frac{
\ell_R(\Tor^R_s(M,\rphir))}{
\ell_R(\Tor^R_s(M,k))} =\ell_R(R) < \infty .
$$
\end{proposition}
\begin{proof}
Let $F_{\bullet} \to M$ be a minimal free resolution of $M$ with 
$F_s=R^{\oplus l_s}$. Then 
$\ell_R(\Tor^R_s(M,k)) = \ell_R(k^{\oplus l_s})=l_s$. On the other
hand, in the complex $\fphir$, all the differentials are zero. 
Hence
$\ell_R(\Tor^R_s(M,\rphir)) = 
\ell_R(F_s)= l_s\ell_R(R)$.
Thus,  
$
\lim_{s \to \infty}\frac{
\ell_R(\Tor^R_s(M,\rphir))}{
\ell_R(\Tor^R_s(M,k))} = \ell_R(R),
$
which is finite since $R$ is Artinian.
\end{proof}

\section{Positive depth and non-rigidity}
 Let $R$ be a local ring of positive depth and 
let $c_\y(R)$ be as defined in~\ref{def1}, where 
$\y = \{y_1, \ldots, y_d\}$ 
is a maximal $R$-regular sequence.
Then the following non-rigidity property holds for
$\Tor ^R_* (-,\rphir)$.
\begin{proposition}\label{prop-not-rigid}
Let $R$ be a local ring with $\depth R = d >0$ and let $M$ be a
finitely generated $R$ module of infinite projective dimension.  If
$\Tor ^R_j (M,\rphir) = 0$ for some $j>0$ and $r > \log _p c_\y (R)$ then
there exists an $n>j$ such that $\Tor ^R_n (M,\rphir)\not = 0$.
\end{proposition}
\begin{proof}
First note that $\y$ is a
 regular sequence on
$\rphir$ since $\{y_1^{p^r}, \ldots, y_d^{p^r}\}$ is again an
$R$-regular sequence.
If $\Tor ^R_{j+1} (M,\rphir) =0$ then from the
long exact sequence 
\begin{equation}\label{eq*}
\cdots \to \Tor ^R_{j+1} (M,\rphir) \to \Tor ^R_{j+1} 
(M,\frac{\rphir}{y_1 \rphir}) \to \Tor ^R_j (M,\rphir) = 0 \to \cdots,
\end{equation}
associated to the short exact sequence 
$$
0 \to \rphir \stackrel{y_1}{\to}\rphir \to \frac{\rphir}{y_1 \rphir} \to 0, 
$$
we conclude that $\Tor ^R_{j+1}(M,\frac{\rphir}{y_1 \rphir}) =0$. 
 Now consider the long exact sequence
\begin{equation}\label{eq**}
\cdots \to \Tor ^R_{j+2} (M,\frac{\rphir}{y_1 \rphir}) \to \Tor ^R_{j+2} 
(M,\frac{\rphir}{(y_1, y_2) \rphir}) \to 
\Tor ^R_{j+1}(M,\frac{\rphir}{y_1 \rphir}) = 0 \to \cdots 
\end{equation}
associated to the short exact sequence 
$$ 
0 \to \frac{\rphir}{y_1
\rphir} \stackrel{y_2}{\to} \frac{\rphir}{y_1 \rphir} \to
\frac{\rphir}{(y_1, y_2) \rphir} \to 0.  
$$ 
From the argument above we get that if $\Tor ^R_{j+2} (M,\rphir)
=0$ then $\Tor ^R_{j+2} (M,\frac{\rphir}{y_1 \rphir}) =0$, and hence
by (\ref{eq**}), $\Tor ^R_{j+2} (M,\frac{\rphir}{(y_1, y_2) \rphir}) =
0$.  Inductively, if $\Tor ^R_{j+i} (M,\rphir) =0$ for $i=1, \ldots,
t$ then $\Tor ^R_{j+t} (M,\frac{\rphir}{(y_1, \ldots, y_t) \rphir}) =0$.

Let $F_\bullet$ be a minimal free resolution of $M$. Then 
$
\Tor ^R_{j+d} (M, ^{\phi^r}\!\!\bar{R}) \cong
H_{j+d} ( ^{\phi^r}\!\! \bar{F}_{\bullet}),  $ 
where $\bar{R} = \frac{R}{\y R}$ and 
$\bar{F}_{\bullet}= F_{\bullet}\otimes_R \bar{R}$. 
By assumption,  $ ^{\phi^r}\!\! \bar{F}_{j+d} \neq 0 $ and
hence, by 
Proposition~\ref{prop-depth0},
 $\Tor ^R_{j+d} (M,^{\phi^r}\!\!\bar{R}) \not =0$.  Thus, by
the inductive argument, we cannot have $\Tor ^R_{j+i}
(M,^{\phi^r}\!\!{R}) = 0$ for all $i=1, \ldots, d$.
\end{proof}
\begin{remark} Note additionally that we have  proved more.  
Namely, that if $\depth (R)>0$ and $M$ is a finitely generated 
$R$-module such that for some $r> \log_p c_\y(R)$ and $n\geq 1$, 
$\Tor ^R_{j} (M,^{\phi^r}\!\!{R}) =0$  for
all $n\leq j \leq n + \depth R$,
 then $M$ has finite projective dimension. As was pointed out in~\cite{AM}, 
this result is implicit in~\cite[(2.6)]{KL}.
It sharpens the results of Peskine and Szpiro \cite[(1.7)]{PS} and
Herzog \cite[(3.1)]{H} which state 
that
the following conditions are equivalent for a finitely generated
$R$-module $M$:
\begin{enumerate}
\item $M$ has finite projective dimension;
\item $\Tor ^R_{j}(M,^{\phi^r}\!\!{R}) = 0$ for all $j, r \geq 1$;
\item $\Tor ^R_{j}(M,^{\phi^r}\!\!{R}) = 0$ for all $j\geq 1$ and 
infinitely many $r$'s.
\end{enumerate}
We note that it may still be possible for $\Tor ^R_{j}(-,^{\phi^r}\!\!{R})$
to be rigid for $r < \log_p c_\y(R)$, or vacuously rigid for 
$r > \log_p c_\y(R)$, that is, $\Tor ^R_{j}(M,^{\phi^r}\!\!{R})$ may
not vanish non-trivially.
\end{remark}
\section{Examples}
In this section we give a few examples of rings that satisfy 
condition~\eqref{e:condition}. Let $k$ be a field of characteristic $p$.

\begin{example}\label{notCM} 
Let $R$ be the quotient ring 
$R=k[[X,Y]]/(XY, X^2)$. Then $R$ is a local one dimensional 
non-Cohen Macaulay  ring with
maximal ideal $m=(x,y)$, where
 $x$ and $y$ are (respectively) the images of $X$ and $Y$  in $R$.
 The ring $R$
satisfies condition~\eqref{e:condition} since $m \subset \ann(x)$
and $x \not\in m^p$. 
\end{example}

\begin{example}
One can make variations on Example~\ref{notCM}. For example, we 
can increase the dimension of the ring by adding variables, say, 
take $R=k[[X,Y,Z]]/(X^2,XY, XZ)$ to get a non-Cohen Macaulay
 ring of dimension 2 with $m \subset \ann(x)$ and $x \not\in m^p$.
\end{example}

\begin{example}(cf.~\cite{Q})
Let $R$ be the graded ring $k \oplus V \oplus k \oplus 0 \oplus 0 \dots$,
where $V$ is a $k$-vector space of dimension at least 3, and the
multiplication on $R$ is defined by some non-degenerated quadratic form
on $V$. Let $m$ be the irrelevant maximal ideal of $R$. 
Then $R$ is an Artinian Gorenstein ring which is not a
complete intersection, with $m^3=0$ and condition~\eqref{e:condition}
is also satisfied for $p=2$ since $m^2V=0$ but $V \not\subset m^2$.
\end{example}


\end{document}